\documentclass{proc-l}

\usepackage{hyperref}

\newtheorem{theorem}{Theorem}[section]
\newtheorem{proposition}[theorem]{Proposition}
\newtheorem{corollary}[theorem]{Corollary}
\newtheorem{lemma}[theorem]{Lemma}

\theoremstyle{definition}

\newtheorem{example}[theorem]{Example}

\newtheorem{problem}[theorem]{Problem}

\theoremstyle{remark}

\numberwithin{equation}{section}

\newcommand{\omegabar}{\omega^{\ast}}

\newcommand{\HB}{\mathcal{HB}}

\DeclareMathOperator{\im}{Im}

\begin{document}

\title[Sums with real zeros]{Sums of entire functions having only real zeros}

\author[S.~R.~Adams]{Steven~R.~Adams}

\author[D.~A.~Cardon]{David~A.~Cardon}
\address{Department of Mathematics, Brigham Young University, Provo, UT 84602, USA}
\email{cardon@math.byu.edu}

\subjclass[2000]{Primary 30C15. Secondary 30D05.}

\date{August 9, 2006}

\commby{Juha M. Heinonen}

\begin{abstract}
We show that certain sums of products of Hermite-Biehler entire
functions have only real zeros, extending results of Cardon. As
applications of this theorem, we construct sums of exponential
functions having only real zeros, we construct polynomials
having zeros only on the unit circle, and we obtain the
three-term recurrence relation for an arbitrary family of real
orthogonal polynomials. We discuss a similarity of this result
with the Lee-Yang circle theorem from statistical mechanics.
Also, we state several open problems.
\end{abstract}

\maketitle

\section{Introduction}

In the study of entire functions, the Hermite-Biehler class of
functions plays a particularly important role. An entire function
$\omega(z)$ is said to belong to the Hermite-Biehler class, denoted
$\HB$, if all of its roots belong to the open upper half plane
$\im(z) > 0$ and if
\begin{equation}\label{eqn:HBdefinition}
\left|\frac{\omega(z)}{\omegabar\!(z)}\right| < 1 \qquad \text{for} \qquad \im(z) > 0
\end{equation}
where $\omegabar(z) = \overline{\omega(\bar{z})}$.

Hermite-Biehler entire functions naturally appear in many areas of
mathematics, and the problem of studying functions with zeros in a
half plane is closely related to studying pairs of entire functions
with interlacing zeros on a line. For an extensive and detailed
discussion of the Hermite-Biehler class, we refer the reader to
chapter VII of~\cite{Levin1980} and the references therein.

A simple, but significant, observation is that the roots of
\begin{equation}\label{eqn:simpleobservation}
\omega(z) + \omegabar(z)
\end{equation}
are necessarily real. This follows from
inequality~\eqref{eqn:HBdefinition}. In this paper, we generalize
this observation by showing, in Theorem~\ref{theorem:main}, that
certain more complicated sums of products of Hermite-Biehler
functions also have only real zeros. This extends results of Cardon
in \cite{Cardon2004} which were used in~\cite{Cardon2005Fourier} to
construct Fourier transforms having only real zeros.

We will now give a precise statement of Theorem~\ref{theorem:main},
our main result. Let $\omega_1(z), \ldots, \omega_n(z) \in \HB$. Let
$G(z)$ be an entire function of genus $0$ or $1$ that is real on the
real axis and has only real zeros. This is equivalent to saying that
$G(z)$ has a product representation of the form $G(z)=cz^q e^{\alpha
z} \prod (1-z/\alpha_m)e^{z/\alpha_m}$ where $q$ is a nonnegative
integer, $c$ and $\alpha$ are real, and the $\alpha_m$ are the
nonzero real zeros of $G$. In the case of either genus $0$ or $1$,
the sum $\sum \alpha_m^{-2}$ is finite. Let $a_1,\ldots,a_n$ be
positive real numbers. Let $T=\{1,2,\ldots,n\}$. For a subset $S$ of
$T$, let $S'$ denote its complement in $T$. Define $H_n(z)$ by
\begin{equation}\label{eqn:Hn}
H_n(z) = \sum_{S \subseteq T} G\Bigl(-\sum_{k \in S'} i\,a_k + \sum_{\ell \in S} i\,a_{\ell}\Bigr)
\prod_{k \in S'} \omegabar_k(z) \prod_{\ell \in S} \omega_{\ell}(z).
\end{equation}
For example,
\begin{align*}
H_2(z) & = G(-ia_1 -ia_2) \omegabar_1(z)  \omegabar_2(z) + G(-ia_1+ia_2)\omegabar_1(z) \omega_2(z) \\
       & \quad +G(ia_1 -ia_2) \omega_1(z) \omegabar_2(z) + G(ia_1+ia_2) \omega_1(z) \omega_2(z).
\end{align*}

The main result of this paper is the following theorem:

\begin{theorem}  \label{theorem:main}
All of the zeros of $H_n(z)$ are real.
\end{theorem}

The organization of the rest of this paper is as follows: In
\S\ref{section:PolyaLeeYang}, we explain some of the motivating
ideas that led us to discover Theorem~\ref{theorem:main}. This
involves combining a result of P\'olya in his study of the Riemann
$\xi$-function from analytic number theory with a result of Lee and
Yang in their study of statistical mechanics. In
\S\ref{section:Proof}, we present the proof of
Theorem~\ref{theorem:main}. Finally, in \S\ref{section:Examples}, we
give examples and applications of this theorem, and we state several
open problems suggested by our studies.

\section{A theorem of Polya and a theorem of Lee and Yang}  \label{section:PolyaLeeYang}

In 1926 P\'olya was attempting to understand the zeros of the
Riemann zeta function. His paper includes the following observation:

\begin{proposition}[P\'olya~\cite{Polya1926}, Hilfssatz II] \label{prop:HilfssatzII}
Let $a>0$ and let $b$ be real. Assume $G(z)$ is an entire function
of genus $0$ or $1$ that is real for real $z$, has at least one real
zero, and has only real zeros. Then the function
\[
G(z-ia)e^{-ib}+G(z+ia) e^{ib}
\]
has only real zeros.
\end{proposition}

Given sequences of real numbers $a_1,a_2,\ldots$ and
$b_1,b_2,\ldots$ it seems natural to iteratively apply the process
of P\'olya's Hilfssatz II (Proposition~\ref{prop:HilfssatzII}) to
form a sequence of functions $F_k(s,z)$ of complex variables $s$ and
$z$ as follows:
\begin{align*}
F_0(s,z) & = G(s), \\
F_k(s,z) & = F_{k-1}(s-ia_k,z) e^{-ib_k z} + F_{k-1}(s+ia_k,z) e^{ib_k z} \quad \text{for $k>0$}.
\end{align*}
For example,
\begin{align*}
F_2(s,z) & = G(s-ia_1-ia_2) e^{iz(-b_1-b_2)} + G(s+ia_1-ia_2)e^{iz(b_1-b_2)} \\
& \qquad + G(s-ia_1+ia_2) e^{iz(-b_1+b_2)} + G(s+ia_1+ia_2)e^{iz(b_1+b_2)}.
\end{align*}

By Proposition~\ref{prop:HilfssatzII}, the functions $F_k(s,z_0)$
for fixed real $z_0$, if nonzero, are of genus $0$ or $1$, have only
real zeros, and are real for real $s$. By taking limits of such
expressions, Cardon~\cite{Cardon2002b} and
Cardon-Nielsen~\cite{Cardon2002a} classified certain distribution
functions $\mu$ such that $\int_{-\infty}^{\infty} G(s-it) d\mu(t)$
has only real zeros. On the other hand, $F_k(s,z)$, as a function of
$z$ for fixed $s$, is a Fourier transform relative to a discrete
measure. In Cardon~\cite{Cardon2004}, it was shown that $F_k(0,z)$
has only real zeros. By taking limits of expressions of the form
$F_k(s,z)$, Cardon~\cite{Cardon2005Fourier} characterized certain
distribution functions $\mu$ such that the Fourier transform
$\int_{-\infty}^{\infty} G(it) e^{izt} d\mu(t)$ has only real zeros.

We may generalize to a greater degree by observing that $\omega(z) =
e^{iz}$ belongs to the Hermite-Biehler class since all of its roots
(there are none) lie in the upper half plane and since
\[
\left|\frac{\omega(z)}{\omegabar(z)}\right| = \left|\frac{e^{iz}}{e^{-iz}}\right|
 = |e^{2iz}|<1 \qquad \text{for $\im(z)>0$}.
\]
Now we iteratively apply P\'olya's Hilffsatz II
(Proposition~\ref{prop:HilfssatzII}) using Hermite-Biehler functions
$\omega_k(z)$, rather than exponential functions, to obtain a
sequence $H_k(s,z)$ of functions of two variables as follows:
\begin{align*}
H_0(s,z) & = G(s), \\
H_n(s,z) & = H_{n-1}(s-ia_n,z) \omegabar_k(z) + H_{n-1}(s+ia_n,z) \omega_k(z), \quad n>0.
\end{align*}
For example,
\begin{align*}
H_2(s,z) & = G(s-ia_1 -ia_2) \omegabar_1(z)  \omegabar_2(z) + G(s-ia_1+ia_2)\omegabar_1(z) \omega_2(z) \\
       & \quad +G(s+ia_1 -ia_2) \omega_1(z) \omegabar_2(z) + G(s+ ia_1+ia_2) \omega_1(z) \omega_2(z).
\end{align*}
We can also define $H_n(s,z)$ in a nonrecursive fashion. Let
$T=\{1,2,\ldots,n\}$. For a subset $S$ of $T$, let $S'$ denote its
complement in $T$. Then $H_n(s,z)$ is the sum
\begin{equation}\label{eqn:Hksz}
H_n(s,z) = \sum_{S \subseteq T} G\Bigl(s-\sum_{k \in S'} i\,a_k + \sum_{\ell \in S} i\,a_{\ell}\Bigr)
\prod_{k \in S'} \omegabar_k(z) \prod_{\ell \in S} \omega_{\ell}(z).
\end{equation}
Theorem~\ref{theorem:main} says that $H_n(z)=H_n(0,z)$ has only real
zeros.

The key to understanding the zeros of $H_n(s,z)$ is to use the
`de-coupling' procedure of Lee and Yang in~\cite{LeeYang1952}. We
will associate a polynomial in the variables $x_1,\ldots,x_n$ with
the function $H_n(s,z)$. Let
\begin{equation} \label{eqn:MyPn}
\begin{split}
P_0(s; x) & = G(s),  \\
P_n(s;x) & = P_{n-1}(s-ia_n;x) + P_{n-1}(s+ia_n; x) x_n \quad \text{for $n>0$},
\end{split}
\end{equation}
where $x=(x_1,\ldots,x_n)$. Note that $x$ is a vector of variables,
the number of which depends on the subscript in the expression
$P_n$. We could also write
\begin{equation} \label{eqn:multivariatepolynomialP}
P_n(s;x) = \sum_{S \subseteq T} G\Bigl(s-\sum_{k \in S'} i\,a_k + \sum_{\ell \in S} i\,a_{\ell}\Bigr)
\prod_{\ell \in S} x_{\ell}.
\end{equation}
The polynomial $P_n(s;x)$ is fundamentally related to $H_n(z)$:
\begin{equation} \label{eqn:relationshipHandP}
P_n\left(s; \frac{\omega_1(z)}{\omegabar_1(z)},\ldots,\frac{\omega_n(z)}{\omegabar_n(z)}\right)
=
\frac{H_n(s,z)}{\omegabar_1(z) \cdots \omegabar_n(z)}.
\end{equation}

The polynomial $P_n(s;x)$ is similar to a polynomial studied by Lee
and Yang. In 1952, they studied a model for phase transitions in
lattice gases with the property that the zeros of the partition
function for the system all lie on the unit circle in the complex
plane. As re-stated in~\cite[p.~108]{Ruelle1999}, the mathematical
result of Lee and Yang~\cite[Appendix II]{LeeYang1952} is

\begin{proposition}[Lee-Yang] \label{prop:LeeYang}
Let $(A_{ij})_{j \neq i}$ be a family of real numbers such that $-1<
A_{ij}<1$, $A_{ij}=A_{ji}$ for $i=1,\ldots,n$; $j=1,\ldots,n$. We
define a polynomial $\mathcal{P}_n$ in $n$ variables by
\begin{equation} \label{eqn:LeeYangPn}
\mathcal{P}_n(x_1,\ldots,x_n) = \sum_{S \subseteq T} \left(\prod_{i \in S} \prod_{j \in S'} A_{ij}\right) \prod_{\ell \in S} x_{\ell}
\end{equation}
where the summation is over all subsets $S$ of $T=\{1,\ldots,n\}$
and where $S'$ is the complement of $S$ in $T$. Then
$\mathcal{P}_n(x_1,\ldots,x_n)=0$ and $|x_1| \geq 1, \ldots,
|x_{n-1}| \geq 1$ imply $|x_n| \leq 1$.
\end{proposition}

Note the similarity of $P_n(s;x_1,\ldots,x_n)$ in~\eqref{eqn:multivariatepolynomialP}
and $\mathcal{P}_n(x_1,\ldots,x_n)$ in~\eqref{eqn:LeeYangPn}. By
setting each $x_k$ equal to $t$, we obtain the Lee-Yang Circle
Theorem from statistical mechanics:

\begin{corollary}[Lee-Yang Circle Theorem]   \label{corollary:LeeYangCircleTheorem}
Let $\mathcal{P}_n(t) = \mathcal{P}_n(t,\ldots,t)$. Then all of the
zeros of $\mathcal{P}_n(t)$ lie on the unit circle.
\end{corollary}

We will now proceed with the proof of Theorem~\ref{theorem:main}.


\section{Proof of Theorem~\ref{theorem:main}}  \label{section:Proof}

In this section we will show that all of the zeros of $H_n(z) \,\,
(=H_n(0,z))$ are real, thus proving Theorem~\ref{theorem:main}.

We begin with a proposition for the multivariate polynomial
$P_n(s;z)$, defined in equation~\eqref{eqn:multivariatepolynomialP}.
This proposition is similar to Proposition~\ref{prop:LeeYang} of Lee
and Yang and is the most important step in the proof of
Theorem~\ref{theorem:main}.

\begin{proposition} \label{prop:multivariatepolynomial}
Let $G(z)$ be an entire function of genus $0$ or $1$ that has only
real zeros, has at least $n \geq 1$ real zeros (counting
multiplicity), and is real for real $z$. Let $A>0$ and
$a_1>0,\ldots,a_n>0$.
\begin{enumerate}
\item[(i)] If $P_n(iA; x_1,\ldots,x_n)=0$ and $|x_1| \geq 1,\ldots, |x_{n-1}| \geq 1$, then $|x_n| < 1$.
\item[(ii)] If $P_n(-iA; x_1,\ldots,x_n)=0$ and $|x_1| \leq 1,\ldots, |x_{n-1}| \leq 1$, then $|x_n| > 1$.
\end{enumerate}
\end{proposition}

\begin{proof}
We will prove part (i) by induction. The proof of part (ii) is
identical except for reversing some inequalities. Let $n=1$, and
suppose that
\[
0 = P_1(iA; x) = G(iA-ia_1)+G(iA+ia_1)x_1.
\]
The hypotheses on $G$ imply $|G(ic)| < |G(id)|$ for real $c,d$ with
$0 \leq |c| < |d|$. Then $|A-a_1|<|A+a_1|$ implies $|x_1| =
\left|\tfrac{G(iA-ia_1)}{G(iA+ia_1)}\right|<1$.

Now assume that $n \geq 2$ and that the theorem holds for $P_k$ with
$1 \leq k < n$. If, by way of contradiction, the theorem is false
for $P_n$, then there exists a solution $x_1,\ldots,x_n$ of the
equation $0 = P_n(iA; x_1,\ldots,x_n)$ such that $|x_k| \geq 1$ for
$1 \leq k \leq n$. We will show this leads to the existence of
another solution $u_1,\ldots,u_n$ such that
\[
|u_1| = 1, \ldots, |u_{n-1}|=1, |u_n| \geq 1.
\]
From the definition of $P_n$,
\begin{align*}
0 & = P_n(iA; x_1,\ldots,x_n) \\
& = P_{n-2}(iA-ia_{n-1}-ia_{n};x) + P_{n-2}(iA-ia_{n-1}+ia_n;x) x_n \\
& \qquad + P_{n-2}(iA+ia_{n-1}-ia_n;x)x_{n-1} + P_{n-2}(iA+ia_{n-1}+ia_{n};x)x_{n-1}x_n.
\end{align*}
By the induction hypothesis,
$P_{n-2}(iA+ia_{n-1}+ia_n;x_1,\ldots,x_{n-2})\neq 0$. The last
equation shows that $x_{n-1}$ and $x_n$ are related through the
fractional linear transformation
\[
x_{n-1} = -\frac{P_{n-2}(iA-ia_{n-1}-ia_n;x) + P_{n-2}(iA-ia_{n-1}+ia_n;x) x_n}{
P_{n-2}(iA+ia_{n-1}-ia_n;x)+P_{n-2}(iA+ia_{n-1}+ia_n;x) x_n}.
\]
As $x_n$ tends to $\infty$, $x_{n-1}$ tends to the value $x'_{n-1}$ where
\[
x'_{n-1} = - \frac{P_{n-2}(iA-ia_{n-1}+ia_n; x)}{P_{n-2}(iA+ia_{n-1}+ia_n;x)}.
\]
For the value $x'_{n-1}$ we have
\begin{align*}
0 & = P_{n-2}(iA-ia_{n-1}+ia_n;x) + P_{n-2}(iA+ia_{n-1}+ia_n;x) x'_{n-1}, \\
0 & = P_{n-1}(iA+ia_n; x_1,\ldots,x_{n-2},x'_{n-1}).
\end{align*}
Since $|x_1|\geq 1,\ldots,|x_{n-2}|\geq 1$, the induction hypothesis
implies that $|x'_{n-1}|<1$. Therefore, by continuity, there exists
a solution $x_1,\ldots,x_{n-2},u_{n-1},\tilde{x}_n$ such that
\[
|x_1| \geq 1, \ldots , |x_{n-2}| \geq 1, |u_{n-1}|=1, |\tilde{x}_n| \geq |x_n| \geq 1.
\]
Repeating this argument with the indices $1,\ldots, n-2$ in place of
$n-1$ results in a solution $u_1,\ldots,u_n$ such that
\[
|u_1|=1,\ldots,|u_{n-1}|=1,|u_n| \geq |x_n| \geq 1.
\]
Then
\begin{align*}
0 & = P_n(iA; u_1,\ldots, u_n), \\
0 & = P_{n-1}(iA-a_n; u_1,\ldots,u_{n-1}) + P_{n-1}(iA+ia_n;u_1,\ldots,u_{n-1})u_n.
\end{align*}
Let $f_k(s)$ for $1 \leq k \leq n-1$ be defined recursively by
\begin{align*}
f_1(s) & = G(s-ia_1)+G(s+ia_1)u_1, \\
f_k(s)& = f_{k-1}(s-ia_k)+f_{k-1}(s+ia_k) u_k.
\end{align*}
Note that
\begin{align*}
f_{n-1}(iA-ia_n) & = P_{n-1}(iA-ia_n; u_1,\ldots,u_{n-1}) \quad \text{and} \\
f_{n-1}(iA+ia_n) & = P_{n-1}(iA+ia_n; u_1,\ldots,u_{n-1}).
\end{align*}

By Lemmas~\ref{lemma:finitelymanyzeros}
and~\ref{lemma:infinitelymanyzeros} (which are stated after this
proof), $f_k(s)$ for $1 \leq k \leq n-1$ is a function of genus $0$
or $1$ that is real for real $s$, has only real zeros, and has at
least one real zero. These conditions, especially the fact that
$f_k(s)$ has at least one real zero, imply that
$|f_{n-1}(iA-ia_n)|<|f_{n-1}(iA+ia_n)|$ since $|A-a_n|<|A+a_n|$.
This gives
\[
|u_n| = \left| \frac{f_{n-1}(iA-ia_n)}{f_{n-1}(iA+ia_n)} \right| <1,
\]
where the inequality is strict. This contradicts the fact that
$|u_n| \geq |x_n| \geq 1$. Therefore, the assumption that $|x_n|
\geq 1$ is false, and we conclude that $|x_n| <1$.
\end{proof}

The following two technical lemmas were used in the proof of
Proposition~\ref{prop:multivariatepolynomial}. They describe the
number of zeros of $G(z-ia)e^{-ib}+G(z+ia)e^{ib}$.

\begin{lemma} \label{lemma:finitelymanyzeros}
Let $G(z)$ be an entire function of genus $0$ or $1$ that is real
for real $z$, has only real zeros, and has at least $n \geq 1$
zeros (counting multiplicity). Let $a>0$ and let $b$ be real. Then
\[G(z-ia)e^{-ib}+G(z+ia)e^{ib}\] is also a function of genus $0$
or $1$ that is real on the real axis, has only real zeros, and has
at least $n-1$ zeros.
\end{lemma}

\begin{proof}
Assume $G$ has $n \geq 1$ zeros. Then $G(z)$ is of the form
\[
G(z) = e^{\alpha z}(c_n z^n+c_{n-1}z^{n-1}+\cdots+c_1z+c_0)
\]
where $c_n \not=0$. Let
\[
H(z) = G(z-ia)e^{-ib}+G(z+ia)e^{ib}.
\]
Then by expanding and collecting powers of $z$, we obtain
\[
H(z) = e^{\alpha z}(d_n z^n + d_{n-1} z^{n-1} + \cdots + d_1 z + d_0)
\]
where
\[
d_n = 2c_n\cos(b+\alpha a) \quad \text{and} \quad
d_{n-1} = 2c_{n-1}\cos(b+\alpha a)-2anc_n \sin(b+\alpha a).
\]
If $\cos(b+\alpha a) \not= 0$, then $d_n \not=0$ and $H(z)$ has
$n$ zeros counting multiplicities. If $\cos(b+\alpha a) = 0$, then
$d_n=0$ but $d_{n-1}\not=0$ in which case $H(z)$ has $n-1$ zeros
counting multiplicities.
\end{proof}

\begin{lemma} \label{lemma:infinitelymanyzeros}
Let $G(z)$ be an entire function of genus $0$ or $1$ that is real
for real $z$, has only real zeros, and has infinitely many zeros.
Let $a>0$ and let $b$ be real. Then
\[G(z-ia)e^{-ib}+G(z+ia)e^{ib}\] is also a function of genus $0$
or $1$ that is real on the real axis, has only real zeros, and has
infinitely many zeros.
\end{lemma}

\begin{proof}
We begin with a simple observation about real entire functions of
genus $0$ or $1$. Let $\phi$ be such a function. Then $\phi$ may
be represented as
\[
\phi(z) = c z^m e^{\alpha z}\prod_{k} (1-z/\alpha_k)e^{z/\alpha_k}
\]
where $c$, $\alpha$, and $\alpha_k$ are real and $m$ is a
nonnegative integer.  For any real $T$,
\[
|\phi(iT)|^2 = c^2 T^{2m} \prod_{k} (1+T^2/\alpha_k^2).
\]
Thus $\phi(z)$ has infinitely many zeros if and only if
$|\phi(iT)|^2$ grows more rapidly than any power of $T$.

Now let
\[
H(z) = G(z-ia)e^{-ib}+G(z+ia)e^{ib}.
\]
By making a change of variable, if necessary, there is no loss of
generality in assuming that $G$ is of the form
\[
G(z) = e^{\alpha z} \prod_{k=1}^{\infty} (1-z/\alpha_k) e^{z/\alpha_k}
\]
where the $\alpha_k$ are the real zeros of $G$. Let
\[
g(z)
=\prod_{k=1}^{\infty} (1-z^2/\alpha_k^2).
\]
For real $T$, $g(iT)=|G(iT)|^2$. Since $g(z)$ satisfies the
conditions of Proposition~\ref{prop:HilfssatzII}, the derivative
$g'(z)=\lim_{h \rightarrow 0} \frac{g(z+ih)-g(z-ih)}{2ih}$ also
satisfies the conditions of Proposition~\ref{prop:HilfssatzII}.
Thus, $g'(z)$ is of the form
\[
g'(z) = z \prod_{k=1}^{\infty} (1-z^2/\beta_k^2),
\]
where the $\beta_k$ are real and $\sum \beta_k^{-2}<\infty$.

By the observation at the beginning of the proof, since $G(z)$ has
infinitely many zeros, $g(iT)=|G(iT)|^2$ grows more rapidly than
any power of $T$. Similarly, for fixed real $a$, both
$g(iT+ia)=|G(iT+ia)|^2$ and $g(iT-ia)=|G(iT-ia)|^2$ grow more
rapidly than any power of $T$.  By showing that the difference
\[
g(iT+ia)-g(iT-ia) = |G(iT+ia)|^2-|G(iT-ia)|^2
\]
grows more rapidly than any power of $T$, we may conclude that
$|H(iT)|^2$ also grows rapidly and hence that $H(z)$ has
infinitely many zeros.

By the mean value theorem of calculus there exists a real number
$a_T$ depending on $T$ in the interval $(-a,a)$ such that
\[
g(iT+ia)-g(iT-ia) = 2ai g'(iT+ia_T) .
\]
Since $g(iT)$ and $ig'(iT)$ are increasing functions of positive
$T$,
\[
g(iT+ia)-g(iT-ia) \geq 2ai g'(iT-ia)
\]
for all $T \geq a$. But the right hand side grows more rapidly
than any power of $T$. Therefore, $|G(iT+ia)|^2$ grows
sufficiently more rapidly than $|G(iT-ia)|^2$ to conclude that
$|H(iT)|^2 =| G(iT+ia)e^{ib}+G(iT-ia)e^{-ib}|^2$ grows more
rapidly than any polynomial as $T$ becomes large. Thus, $H(z)$ has
infinitely many zeros.
\end{proof}

We may now apply the information in
Proposition~\ref{prop:multivariatepolynomial} about the multivariate
polynomial $P_n(s;x_1,\ldots,x_n)$ to the function $H_n(s,z)$,
defined in equation~\eqref{eqn:Hksz}. This results in a lemma which
is Theorem~\ref{theorem:main}, except for the requirement that
$G(z)$ have $n$ real zeros.

\begin{lemma} \label{lemma:almostthemaintheorem}
Let $G(z)$ be an entire function of genus $0$ or $1$ that has only
real zeros, has at least $n \geq 1$ real zeros (counting
multiplicity), and is real for real $z$. Suppose
$A>0,a_1>0,\ldots,a_n>0$ and let $\omega_1(z),\ldots,\omega_n(z) \in
\HB$.
\begin{enumerate}
\item[(i)]
If $H_n(iA,z) = 0$, then $\im(z)>0$.
\item[(ii)]
If $H_n(-iA,z)=0$, then $\im(z)<0$.
\item[(iii)]
If $H_n(0,z)=0$, then $z$ is real.
\end{enumerate}
\end{lemma}

\begin{proof}
From the definition of $H_n(s,z)$, it is immediate that $H_n(iA, z)
= 0$ if and only if $H_n(-iA,\bar{z})=0$. Hence, (i) and (ii) are
equivalent. Recall from equation~\eqref{eqn:relationshipHandP} that
$P_n$ and $H_n$ are related by
\[
P_n\left(s; \frac{\omega_1(z)}{\omegabar_1(z)},\ldots,\frac{\omega_n(z)}{\omegabar_n(z)}\right)
=
\frac{H_n(s,z)}{\omegabar_1(z) \cdots \omegabar_n(z)}.
\]
If $\im(z) \geq 0$, then $|\omega_k(z)/\omegabar_k(z)| \leq 1$ for
$1 \leq k \leq n$. Thus,
Proposition~\ref{prop:multivariatepolynomial}(ii) implies $H(-iA,z)
\not=0$. So, $H_n(-iA,z)=0$ implies $\im(z)<0$. This proves (ii) and
consequently (i).

By part (i), since $\lim_{A\rightarrow 0^+} H_n(iA,z)=H_n(0,z)$ is
uniform on compact sets, the zeros of $H_n(0,z)$ satisfy $\im(z)
\geq 0$. Similarly, by part (ii), the zeros of $\lim_{A \rightarrow
0^+} H_n(-iA,z)=H_n(0,z)$ satisfy $\im(z) \leq 0$. Therefore, the
zeros of $H_n(0,z)$ are purely real. This proves (iii).
\end{proof}

The final step of the proof of Theorem~\ref{theorem:main} is to
improve Lemma~\ref{lemma:almostthemaintheorem} by removing the
restriction that $G(z)$ must have at least $n \geq 1$ real zeros.
The case not covered is when $G(z)$ is of the form
\[
G(z) = c z^q e^{\alpha z} \prod_{m=1}^k(1-z/\alpha_m)
\]
where $q$ is a nonnegative integer, $c$ and $\alpha$ are real, the
$\alpha_m$ are the nonzero real zeros of $G$, and $0 \leq k < n$.
For positive $N$ let
\[
G_N(z) = (1-z/N)^{n-k} \cdot c z^q e^{\alpha z} \prod_{m=1}^k(1-z/\alpha_m).
\]
Then $G_N(z)$ has $n$ real roots. Let $H_{N,n}(z)$ be the sum of
Hermite-Biehler functions corresponding to $G_N(z)$:
\[
H_{N,n}(z) = \sum_{S \subset T} G_N\Bigl( - \sum_{k \in S'} i a_k + \sum_{\ell \in S} i a_{\ell}\Bigr)
\prod_{k \in S'} \omegabar_k(z) \prod_{\ell \in S} \omega_{\ell}(z).
\]
By Lemma~\ref{lemma:almostthemaintheorem}, $H_{N,n}(z)$ has only
real zeros. Since $\lim_{N \rightarrow \infty} H_{N,n}(z) = H_n(z)$
and the limit is uniform on compact sets, the zeros of $H_n(z)$ are
real.       This completes the proof of Theorem~\ref{theorem:main}.


\section{Applications, Examples, and Questions} \label{section:Examples}
We conclude the paper by providing several examples and applications
of Theorem~\ref{theorem:main}. We also state several open problems.

\begin{example}[Cardon~\cite{Cardon2004}, Theorem 1]  \label{theorem:exponentialsums}
In Theorem~\ref{theorem:main}, let $\omega_k(z) = e^{i b_k z}$. Then
we obtain the following theorem: Let $a_1,a_2,a_3,\ldots$ and
$b_1,b_2,b_3,\ldots$ be positive. Let $G$ be an entire function of
genus $0$ or $1$ that is real on the real axis and has only real
zeros. Let $H_n(z)$ be defined by
\[
H_n(z) = \sum G( \pm ia_1 \pm i a_2 \pm \cdots \pm ia_n) e^{iz(\pm b_1 \pm \cdots \pm b_n)}
\]
where the summation is over all $2^n$ possible plus and minus sign
combinations, the same sign combination being used in both the
argument of $G$ and in the exponent. Then $H_n(z)$ has only real zeros.
\end{example}

\begin{example}
By taking limits of the sums of exponential functions in the
previous example, it is possible to classify certain distribution
functions $\mu(t)$ such that the Fourier transform
\[
H(z) = \int_{-\infty}^{\infty} G(it) e^{izt} d\mu(t)
\]
has only real zeros. This is explained in detail in
Cardon~\cite{Cardon2005Fourier}.
\end{example}

\begin{example}[Circle Theorem]  \label{example:MyCircleTheorem}
In the polynomial $P_n(s;x_1,\ldots,x_n)$ defined in
equation~\eqref{eqn:multivariatepolynomialP}, let $s=0$ and let each
$x_k$ equal $t$. We obtain a polynomial $P_n(t) = P_n(0;t,\ldots,t)$
of a single variable. This polynomial can be described by
\[
P_n(t) = \sum_\sigma G( \sigma \cdot(ia_1,\ldots,ia_n)) t^{|\sigma|}
\]
where the summation is over all $2^n$ vectors $\sigma$ of the form
$(\pm 1, \ldots, \pm 1)$, $|\sigma|$ represents the number of plus
signs in the vectors $\sigma$, and $\sigma \cdot (ia_1,\ldots,ia_n)$
is the ordinary dot product. For example,
\[
P_2(t) = G(-ia_1-ia_2) + G(-i a_1 + i a_2)t + G(i a_1 - i a_2)t +
G(i a_1 + i a_2) t^2.
\]
By Proposition~\ref{prop:multivariatepolynomial}, all the zeros of
$P_n(t)$ lie on the unit circle in the complex plane. This fact
resembles the Lee-Yang Circle Theorem
(Corollary~\ref{corollary:LeeYangCircleTheorem}).
\end{example}

\begin{example}[Orthogonal polynomials]   \label{example:orthogonal}
Let $p_0(z),p_1(z),p_2(z),\ldots$ be any family of real orthogonal
polynomials with $\deg p_n(z) = n$ normalized so that the leading
coefficient of each $p_n(z)$ is positive. The well known three-term
recurrence relation states that there exist real constants
$A_n,B_n,C_n$ such that
\begin{equation} \label{eqn:threeterm}
p_n(z) = (A_n z + B_n) p_{n-1}(z) - C_n p_{n-2}(z), \qquad n=2, 3, 4, \ldots
\end{equation}
where $A_n>0$ and $C_n>0$ (Szeg\"o~\cite[Thm. 3.2.1]{Szego}). We
will show that Theorem~\ref{theorem:main} reproduces
equation~\eqref{eqn:threeterm}. It is known that the zeros of
$p_{n-2}(z)$ and $p_{n-1}(z)$ are real and interlace and that
$p_{n-2}(z)p'_{n-1}(z) - p'_{n-2}(z)p_{n-1}(z)>0$ for real $z$ (see
\cite{Szego}, \S3.3). On the other hand, the Hermite-Biehler theorem
for polynomials states that if $p(z)$ and $q(z)$ are real
polynomials with real interlacing roots and if
$p(z)q'(z)-p'(z)q(z)>0$ for some real $z$, then $p(z)+iq(z) \in \HB$
(Levin~\cite[p.~305]{Levin1980}). Consequently, the polynomial
\[
\omega_1(z) = p_{n-2}(z) + i p_{n-1}(z)
\]
belongs to $\HB$. Let $\omega_2(z)=z-(-B_n/A_n+i)$. Note that
$\omega_2(z) \in \HB$. Let $G(z)=-z$, $a_1=A_n/4$, and $a_2=C_n/4$.
A short calculation gives
\begin{align*}
H_2(z)
& =
G(-ia_1 -ia_2) \omegabar_1(z)  \omegabar_2(z) + G(-ia_1+ia_2)\omegabar_1(z) \omega_2(z) \\
 & \qquad +G(ia_1 -ia_2) \omega_1(z) \omegabar_2(z) + G(ia_1+ia_2) \omega_1(z) \omega_2(z) \\
& =
(A_n z + B_n) p_{n-1}(z) - C_n p_{n-2}(z).
\end{align*}
The right hand side of the last expression matches the right hand
side of equation~\eqref{eqn:threeterm}, and by
Theorem~\ref{theorem:main}, it has only real zeros. Unfortunately,
Theorem~\ref{theorem:main} does not imply that the zeros of
$H_2(z)=(A_n z + B_n) p_{n-1}(z) - C_n p_{n-2}(z)$ are simple or
interlace with the zeros of $p_{n-1}(z)$. Hence,
Theorem~\ref{theorem:main} cannot be directly used in an induction
argument to show that the polynomials defined by the
recurrence~\eqref{eqn:threeterm} have only simple real zeros.
\end{example}

We conclude with several questions for further study:
\begin{problem}
The results for the multivariate polynomials
$\mathcal{P}_n(x_1,\ldots,x_n)$ in Proposition~\ref{prop:LeeYang}
and $P_n(s;x_1,\ldots,x_n)$ in
Proposition~\ref{prop:multivariatepolynomial} are remarkably
similar. Is there a way to view these two results as special cases
of a more general theorem? Does either one of these propositions
imply the other?
\end{problem}

\begin{problem}
Is there a physical interpretation for the polynomial $P_n(t)$ in
Example~\ref{example:MyCircleTheorem} similar to the physical
interpretation of the polynomial $\mathcal{P}_n(t)$ in the Lee-Yang
Circle Theorem (Corollary~\ref{corollary:LeeYangCircleTheorem})?
\end{problem}

\begin{problem}
Determine conditions ensuring that $H_n(z)$ in
Theorem~\ref{theorem:main} has only simple real zeros.
\end{problem}

\begin{problem}
Extend Theorem~\ref{theorem:main} so that it gives information about
the relationship of the zeros of $H_n(z)$ and those of $p_k(z)$ and
$q_k(z)$ where $p_k(z)$ and $q_k(z)$ are real entire functions such
that $\omega_k(z) = p_k(z) + i q_k(z)$ for $1 \leq k \leq n$. In
particular, improve Theorem~\ref{theorem:main}, so that it can be
used to give a direct proof of the simplicity and interlacing of the
zeros of the real orthogonal polynomials discussed in
Example~\ref{example:orthogonal}.
\end{problem}


\bibliographystyle{amsplain}

\end{document}